\documentclass{amsart}
\usepackage{ amsmath, amsthm, amsfonts, hyperref, graphicx, ifpdf}
\usepackage[dvipsnames,usenames]{color}
\hypersetup{
   unicode=false,          
   pdftoolbar=true,        
   pdfmenubar=true,        
   pdffitwindow=false,     
   pdfstartview={},    
   pdftitle={},    
   pdfauthor={},     
   pdfsubject={},   
   pdfcreator={},   
   pdfproducer={}, 
   pdfkeywords={}, 
   pdfnewwindow=true,      
   colorlinks=true,       
   linkcolor=blue,          
   citecolor=blue,        
   filecolor=blue,      
   urlcolor=blue          
}

\newcommand{\R}{{\mathbb R}}

\newtheorem{theorem}{Theorem}[section]
\newtheorem{lemma}[theorem]{Lemma}
\newtheorem{pro}[theorem]{Proposition}
\newtheorem{cor}[theorem]{Corollary}
\newtheorem{prop}[theorem]{Proposition}

\theoremstyle{definition}

\theoremstyle{remark}
\newtheorem{remark}[theorem]{Remark}

\numberwithin{equation}{section}

\newcommand{\ee}{evolution equation}
\newcommand{\mc}{mean curvature}
\newcommand{\mcf}{mean curvature flow}
\newcommand{\maxp}{maximum principle}

\newcommand{\ub}{uniform bound}

\newcommand{\be}{\begin{equation}}
\newcommand{\ene}{\end{equation}}
\newcommand{\br}{\begin{remark}}
\newcommand{\er}{\end{remark}}
\newcommand{\bl}{\begin{lemma}}
\newcommand{\el}{\end{lemma}}
\newcommand{\bcor}{\begin{cor}}
\newcommand{\ecor}{\end{cor}}
\newcommand{\ben}{\begin{enumerate}}
\newcommand{\een}{\end{enumerate}}
\newcommand{\bp}{\begin{proof}}
\newcommand{\ep}{\end{proof}}
\newcommand{\bpo}{\begin{pro}}
\newcommand{\epo}{\end{pro}}
\newcommand{\beq}{\begin{equation*}}
\newcommand{\eeq}{\end{equation*}}
\newcommand{\bear}{\begin{eqnarray*}}
\newcommand{\eear}{\end{eqnarray*}}
\newcommand{\bt}{\begin{theorem}}
\newcommand{\et}{\end{theorem}}

\DeclareMathAlphabet{\mathcal}{OMS}{cmsy}{m}{n}

\numberwithin{equation}{section}

\allowdisplaybreaks

\def\XXint#1#2#3{{\setbox0=\hbox{$#1{#2#3}{\int}$}
    \vcenter{\hbox{$#2#3$}}\kern-.5\wd0}}

\makeatletter
\def\@citestyle{\m@th\upshape\mdseries}
\def\citeform#1{{\bfseries#1}}
\def\@cite#1#2{{%
  \@citestyle[\citeform{#1}\if@tempswa, #2\fi]}}
\@ifundefined{cite }{%
  \expandafter\leqt\csname cite \endcsname\cite
  \edef\cite{\@nx\protect\@xp\@nx\csname cite \endcsname}%
}{}
\makeatother

\begin{document}
\newcommand{\osc}{{\text{osc}}}
\newcommand{\Vol}{{\text{Vol}}}
\newcommand{\V}{{\text{V}}}
\newcommand{\nn}{{\text{n}}}
\newcommand{\cM}{{\cal{M}}}
\newcommand{\Ric}{{\text{Ric}}}
\newcommand{\RE}{{\text{Re }}}
\newcommand{\LL}{{\cal{L}}}
\newcommand{\diam}{{\text {diam}}}
\newcommand{\dist}{{\text {dist}}}
\newcommand{\Area}{{\text {Area}}}
\newcommand{\Length}{{\text {Length}}}
\newcommand{\Energy}{{\text {Energy}}}
\newcommand{\SSS}{{\bold S}}
\newcommand{\K}{{\text{K}}}
\newcommand{\Hess}{{\text {Hess}}}
\def\RR{{\bold  R}}
\def\SS{{\bold  S}}
\def\TT{{\bold  T}}
\def\CC{{\bold C }}
\newcommand{\dv}{{\text {div}}}
\newcommand{\kg}{{\text {k}}}
\newcommand{\expp}{{\text {exp}}}
\newcommand{\e}{{\text {e}}}
\newcommand{\eqr}[1]{(\ref{#1})}
\newcommand{\ec}{\varepsilon_c}
\newcommand{\rc}{\rho_c}
\newcommand{\sztp}{\Sigma^{0,2\pi}}

\parskip1ex






\title[Blow-up of the mean curvature at the first singular time of MCF]{Blow-up of the mean curvature at the first singular time of the mean curvature flow}

\author{Longzhi Lin}
\address[L.~Lin]{Department of Mathematics\\Rutgers University\\110 Frelinghuysen Road\\Piscataway, NJ 08854-8019\\USA}
\email{lzlin@math.rutgers.edu}

\author{Natasa Sesum}
\address[N.~Sesum]{Department of Mathematics\\Rutgers University\\110 Frelinghuysen Road\\Piscataway, NJ 08854-8019\\USA}
\email{natasas@math.rutgers.edu}

\date{}


\begin{abstract}
It is conjectured that the mean curvature blows up at the first singular time of the mean curvature flow in Euclidean space, at least in dimensions less or equal to 7. We show that the mean curvature blows up at the singularities of the mean curvature flow starting from an immersed closed hypersurface with small $L^2$-norm of the traceless second fundamental form (observe that the initial hypersurface is not necessarily convex). As a consequence of the proof of this result we also obtain the dynamic stability of a sphere along the mean curvature flow with respect to the $L^2$-norm.
\end{abstract}

\maketitle


\section{Introduction}
Let $F_0: M^n \to \mathbb{R}^{n+1}$ be an immersion of a closed hypersurface $M$ and evolve it by the mean curvature flow, that is,
\begin{equation}
\label{eq-mcf}
\frac{\partial F}{\partial t} = -H\nu, \qquad F(\cdot,0) = F_0(\cdot)\,.
\end{equation}
Here $H$ is the mean curvature and $\nu$ is the outward unit normal vector of the surface $M_t = F(\cdot, t)$. We denote
by $A = \{a_{ij}\}$ the second fundamental form of $M_t$ and its traceless part
$\text{\AA} = A - \frac{H}{n}g$, whose norm square is given by $$|\text{\AA}|^2 = |A|^2 - \frac{1}{n}H^2 = \frac{1}{n} \sum_{i<j}^n(\kappa_i - \kappa_j)^2$$ where $\kappa_i$'s are the principle curvatures of $M_t$. For an immersed two dimensional surface $\Sigma$ in $\mathbb{R}^n$ the quantity $\int_{\Sigma} |\text{\AA}|^2 d\mu$ is usually referred as Willmore energy of $\Sigma$. The condition $\int_{\Sigma} |\text{\AA}|^2 d\mu < 8\pi$ implies that $\Sigma$ is topologically a sphere (see e.g. \cite{DLM05}).

The traceless second fundamental form measures the roundness of a hypersurface, the smaller it is in a considered norm the closer we are to a sphere in that norm. More precisely, if the $\int_{\Sigma} |\text{\AA}|^2\, d\mu < \epsilon$, we say the hypersurface $\Sigma$ is $\epsilon$-close to the sphere in the $L^2$ sense.
A classical result in differential geometry is a version of Codazzi's theorem which states that every closed, connected and immersed two dimensional surface $\Sigma \in \mathbb{R}^n$ with vanishing traceless second fundamental form, i.e., $\text{\AA} = 0$, is isometric to a round sphere. In fact, it is well-known by now that closed immersed umbilic (i.e., $\text{\AA} = 0$) hypersurfaces in $\mathbb{R}^{n}$ are spheres, see e.g. \cite{Ger}. Note that for surfaces in $\mathbb{R}^3$, De Lellis and M{\"u}ller in \cite{DLM05} generalized Codazzi's theorem by showing the quantitative estimate
$$\inf_{\lambda\in \mathbf{R}} \| A - \lambda \,\text{Id}\|_{L^2(\Sigma)} \leq C \| \text{\AA}\|_{L^2(\Sigma)}
$$for some universal constant $C$. This has been recently generalized to higher dimensional convex hypersurfaces in $\mathbb{R}^n$ in \cite{Per11}. It is not known if this quantitative rigidity still holds for non-convex hypersurfaces. It is in this sense that we say a hypersurface with small $\int |\text{\AA}|^2 d\mu$ is close to round sphere. Geometric flows starting from hypersurfaces with small $\int |\text{\AA}|^2 d\mu$ have been studied a lot. In \cite{HY96} Huisken and Yau used the volume preserving mean curvature flow starting from large spheres (which implies small $|\text{\AA}|$) to construct CMC foliation in the exterior region for asymptotically flat three-manifolds. Kuwert and Sch$\ddot{\text{a}}$tzle \cite{KS01} showed that the Willmore flow (i.e., the gradient flow of the Willmore energy $\int |\text{\AA}|^2 d\mu$) starting from a surface with small $\int |\text{\AA}|^2 d\mu$ in $\mathbb{R}^3$ exists for all time and converges to a round sphere.  The surface area preserving mean curvature flow starting from hypersurfaces with small $\int |\text{\AA}|^2 d\mu$ in $\mathbb{R}^n$ was recently investigated by Huang and the first author in \cite{HL02}.

By the avoidance principle in the mean curvature flow it follows that the flow starting from any closed hypersurface in $\mathbb{R}^{n+1}$ develops singularity at a finite time. In \cite{Hui84} Huisken proved that the norm of the second fundamental form $|A|$ has to blow up at the first singular time.  In \cite{Coo11} Cooper proves the quantity $|A||H|$ needs to blow up at the singularity. It is still an open question whether the mean curvature needs to blow up at the first singular time as well. One expects this to be true at least in dimensions $n \le 7$. In this paper we give a partial answer to this question, namely we prove the following theorem.

\bt\label{MainThm}
Let $M_t^n \subset \mathbb{R}^{n+1}, n\geq 2$ be a smooth compact solution to the {\mcf} \eqref{eq-mcf} for $t\in [0,T)$ with $T<\infty$. Assume there exists a $c_0$ so that $\sup_{M_t} |H|(\cdot,t) \le c_0$ for all $t\in [0,T)$. There exists an $\epsilon>0$ depending only on $n, c_0$, the area of $M_0$, $\max_{M_0}|A|$ and the bound on $\int_{M_0}|\nabla^m A|^2 d\mu$ (for all $m\in [1,\hat{m}]$ for some fixed $\hat{m}\gg1$) such that if
$$\int_{M_0} |\text{\AA}|^2 \,d\mu < \epsilon,$$
then the flow can be smoothly extended past time $T$.
\et

\br
Without the smallness assumption on $\int_{M_0} |\text{\AA}|^2 \,d\mu$, Theorem \ref{MainThm} was proved by Le and the second author in \cite{LS10} for type I singularities. Theorem \ref{MainThm} says that the mean curvature blows up at the singularity of the mean curvature flow if the initial $L^2$-norm of the traceless second fundamental form is small, regardless the type of singularities (but as we shall see, by the following dynamic stability result in Theorem \ref{thm-conv}, the initial smallness of $\int_{M_0}|\text{\AA}|^2 d\mu$ forces the singularities to be of type I).
\er

In \cite{Hui84} Huisken showed that if the initial hypersurface $M\subset \mathbb{R}^{n+1}$ is closed and convex then the normalized mean curvature flow (see equation \eqref{5}) exists forever and exponentially converges to a sphere. This result in particular implies that the normalized mean curvature flow starting at any sufficiently small $C^2$ perturbation of the sphere exponentially converges, as time converges to infinity, to the sphere (since any small $C^2$-perturbation of the sphere is a convex hypersurface). This means the sphere is $C^2$ dynamically stable along the flow, which is in connection to Colding-Minicozzi's classification on $\mathcal{F}$-stable self-shrinkers, see \cite{CM12}. The question is whether there are other, possibly non-convex perturbations of the sphere which exhibit the same type of behavior. For example, in \cite{KS12} it is shown that if $M_0$ is close to an Euclidean $n$-sphere in the Sobolev norm $H^{s}, s> \frac{n}{2}+1$, then the flow contracts to a round point in finite time. Having in mind that the traceless second fundamental form measures the roundness of a hypersurface (in the aforementioned $L^2$ sense), that is, the proximity to a sphere, we prove the following result.

\bt
\label{thm-conv}
Let $M_t^n \subset \mathbb{R}^{n+1}, n\geq 2$ be a smooth compact solution to the {\mcf} \eqref{eq-mcf} for $t\in [0,T)$
with $T < \infty$. Then there exists an
$\epsilon$ depending only on $n$, the area of $M_0$, $\max_{M_0}|A|$ and the bound on $\int_{M_0}|\nabla^m A|^2 d\mu$ (for all $m\in [1,\hat{m}]$ for some fixed $\hat{m}\gg1$) such that if
\be \label{7}
\int_{M_0}|\text{\AA}|^2 \,d\mu < \epsilon,
\ene
then the normalized mean curvature flow \eqref{5} exists for all time and converges exponentially to a round sphere.
\et

The organization of the paper is as follows. In Section \ref{sec-prel} we give necessary preliminaries needed for the proofs of main results. In Section \ref{sec-extension} we prove Theorem \ref{MainThm} and in section \ref{sec-conv} we show Theorem \ref{thm-conv}.

\subsection*{Acknowledgements} The research of the second author is partially supported by NSF grant 1056387. The authors would like to thank Zheng Huang and Nam Le for many helpful discussions.



\section{preliminaries}
\label{sec-prel}

We collect some necessary preliminary results in this section, in particular some well-known evolution equations of several geometric quantities which will be used later.

\bcor (\cite{Hui84})\label{Evol-A-H}
We have the evolution equations for $H$ and $|A|^2$:
\begin{itemize}
\item[(i)] $\frac{\partial}{\partial t} H = \Delta H + |A|^2 H$;
\item[(ii)] $\frac{\partial}{\partial t} |A|^2 = \Delta |A|^2 - 2|\nabla A|^2 + 2|A|^4$\,.
\end{itemize}
Therefore we also have
\begin{itemize}
\item[(iii)] $\frac{\partial}{\partial t} |\text{\AA}|^2 =
\Delta |\text{\AA}|^2 - 2|\nabla \text{\AA}|^2 +
2|A|^2|\text{\AA}|^2$\,,
where $|\nabla \text{\AA}|^2 = |\nabla A|^2 - \frac{1}{n}|\nabla H|^2$\,.
\end{itemize}
\ecor
Immediate corollary of above evolution equations are evolutions of higher order derivatives of the second fundamental form.
\bcor(\cite{Hui84}) \label{hdA-1}
We have the evolution equation for $|\nabla^m A|^2$:
\be \label{hdA-2}
\frac{\partial}{\partial t} |\nabla^m A|^2 = \Delta |\nabla^m A|^2 - 2|\nabla^{m+1} A|^2 + \sum_{i+j+k=m} \nabla^i A\ast \nabla^j A \ast \nabla^k A \ast \nabla^m A \,,
\ene
where $S\ast \Omega$ denotes any linear combination of tensors formed by contraction on $S$
and $\Omega$ by the metric $g$.
\ecor
\noindent Using that
$$\frac{d}{dt}\, d\mu = -H^2\, d\mu,$$
and integrating \eqref{hdA-2} by parts (see \cite{Hui84})  we obtain that
for any $m\geq 1$ we have the estimate
\begin{align}\label{hdA-3}
&\frac{d}{dt}\int_{M_t} |\nabla^m A|^2 \,d\mu + \int_{M_t} |\nabla^m A|^2H^2 \,d\mu\\
\leq &-2\int_{M_t} |\nabla^{m+1} A|^2\,d\mu + C(n,m)\,\max_{M_t}|A|^2 \int_{M_t} |\nabla^m A|^2 \,d\mu\,. \nonumber
\end{align}

Let us also recall the interpolation inequalities for tensors proved by Hamilton.

\bt(\cite{Ham82})\label{Hamiltonlemma2}
Let $M$ be an $n$-dimensional compact Riemannian manifold and $\Omega$ be any tensor on $M$.
\begin{enumerate}
\item[(i)]
Suppose $\frac{1}{p} + \frac{1}{q} = \frac{1}{r}$ with  $r\ge 1$. Then
\beq
\left(\int_{M} |\nabla \Omega|^{2r}\, d\mu\right)^{1/r} \leq (2r-2+n)\,
\left(\int_{M} |\nabla^2 \Omega|^p\, d\mu\right)^{1/p}\left(\int_{M} |\Omega|^{q}\, d\mu\right)^{1/q}\,.
\eeq
\item[(ii)]
If $1\leq i\leq n-1$ and $j \ge 0$  there exists a constant $C = C(n,j)$ which is independent of the metric and connection on $M$ such that
\beq
\int_{M} |\nabla^i \Omega|^{2j/i}\, d\mu \leq C\, \max_{M} |\Omega|^{2(j/i-1)}\int_{M}|\nabla^j \Omega|^2\, d\mu\,.
\eeq
\end{enumerate}
\et

We next state the original Michael-Simon's inequality. A variant of this inequality will be stated in Lemma \ref{MSineq2}.
\bl \label{MSineq1} (\cite{MS73}) Let $M$ be a compact $n$-dimensional hypersurface without boundary, which is smoothly embedded in $\mathbb{R}^{n+1}$. For any Lipschitz function $v\geq 0$ on $M$ we have
\be
\label{eq-MS}
\left( \int_M v^{\frac{n}{n-1}}\, d\mu \right)^{\frac{n-1}{n}} \le C(n)\, \left( \int_M |\nabla v|\, d\mu + \int_M |H| v\, d\mu \right).
\ene
\el

The following result of Topping about the bound on the diameter of an evolving hypersurface will be useful in the proof of the main theorem.
\bl(\cite{Top08})\label{top}
Let $M$ be an $n$-dimensional closed, connected manifold smoothly immersed in
$\mathbb{R}^N$, where $N \geq n+1$. Then the intrinsic diameter and the {\mc} $H$ of $M$ are related by
\beq
\text{diam}\,(M) \leq C(n)\int_M |H|^{n-1}\,d\mu\,.
\eeq
\el

Finally we state the version of the {\maxp} that we will use below (especially in the proof of  Proposition \ref{key1}).
\bt\label{maxPrin} (Maximum principle, see e.g. \cite[Lemma 2.12]{CLN06}) Suppose $u: M \times [0,T] \to \R$ satisfies
\beq
\frac{\partial}{\partial t} u  \leq a^{ij}(t)\nabla_i\nabla_j u + \langle B(t), \nabla u\rangle + F(u)\,,
\eeq
where the coefficient matrix $\left(a^{ij}(t)\right)>0$ for all $t\in [0,T]$, $B(t)$ is a time-dependent vector field and $F$ is a
Lipschitz function. If $u\leq c$ at $t=0$ for some $c>0$, then $u(x,t)\leq U(t)$ for all $(x,t)\in M_t, t\ge 0$, where $U(t)$
is the solution to the following initial value problem:
\beq
\frac{d}{dt} U(t) = F(U) \quad \text{with} \quad U(0) = c\,.
\eeq
\et


\section{Small traceless second fundamental form}
\label{sec-extension}

In this section we prove Theorem \ref{MainThm}.  Our strategy is as follows. We start with a hypersurface with small $L^2$-norm of the traceless second fundamental form. If initially the $L^2$-norms of the second fundamental form $A$ and its derivatives $\nabla^m A$ for $m \in [1, \hat{m}]$ (where $\hat{m} \gg 1$) are bounded by some $\Lambda_0 \gg 1$, using the smallness of the $L^2$-norm of the traceless second fundamental form we show that $|\text{\AA}|(\cdot,t)$ stays uniformly small for short time $t\in (0,T_1]$, where $T_1 = T_1(\Lambda_0) \in (0,1)$. This will imply that (choosing $\Lambda_0 \gg c_0$ where $c_0$ is the uniform bound of the mean curvature)
$$\max\left\{\max_{M_t} |A|, \, \int_{M_t} |\nabla^m A|^2\, d\mu \right\} \le \frac{\Lambda_0}{2},$$
for all $t\in (0,T_1]$ and $m \in [1,\hat{m}]$. Using the uniform bound on the mean curvature and the pointwise smallness of $|\text{\AA}|(\cdot,t)$ for $t\in (0,T_1]$ we also show that $\int_{M_t} |\text{\AA}|^2\, d\mu$ decreases along the flow, so is therefore even smaller than initially. We use that to iterate our arguments starting now at $t = T_1$ instead of $t = 0$
on a time interval of uniform size $T_1>0$. This means that after finitely many iterations we reach time $T$ showing that $|A|$ can not blow up at time $T$ unless the mean curvature does.

In order to carry out our proof we  need the following version of  Michael-Simon's inequality.
\bl \label{MSineq2} Let $M$ be a closed $n$-dimensional hypersurface, smoothly immersed in $\mathbb{R}^{n+1}$. Let $v \geq 0$ be any Lipschitz function on $M$.
We have:
\begin{enumerate}
\item[(i)]
For any $n > 2$,
\be
\label{eq-MS-1}
\left( \int_M v^{\frac{2n}{n-2}}\, d\mu \right)^{\frac{n-2}{n}} \leq C(n)\, \left( \int_M |\nabla v|^2\, d\mu + \int_M H^2 v^2\, d\mu \right)\,.
\ene
\item[(ii)]
For $n = 2$,
\be
\label{eq-MS-2}
\int_M v^2 \le C(n)\, \left(\int_M |\nabla v|^2\, d\mu + \int_M H^2 v^2\, d\mu\right)\,.
\ene
\end{enumerate}
\el
\bp
\begin{enumerate}
\item[(i)]
Apply
 Michael-Simon's inequality \eqref{eq-MS} to  function $w = v^{\frac{2(n-1)}{n-2}}$ to get
$$
\left( \int_M v^{\frac{2n}{n-2}} \, d\mu\right)^{\frac{n-1}{n}} \leq C(n) \left( \int_M |\nabla v| v^{\frac{n}{n-2}} \, d\mu + \int_M |H| v^{\frac{2(n-1)}{n-2}} \, d\mu\right)\,.
$$
It follows by H$\ddot{\text{o}}$lder's inequality that
\begin{align*}
\left( \int_M v^{\frac{2n}{n-2}} \, d\mu\right)^{\frac{n-2}{n}} \leq &\, C(n) \left( \int_M |\nabla v| v^{\frac{n}{n-2}} \, d\mu + \int_M |H| v\cdot v^{\frac{n}{n-2}} \, d\mu\right)^{\frac{n-2}{n-1}} \\
\leq  & \, C(n) \left( \int_M |\nabla v|^2 \, d\mu \right)^{\frac{n-2}{2(n-1)}}\left( \int_M v^{\frac{2n}{n-2}} \, d\mu \right)^{\frac{n-2}{2(n-1)}} \\
& + C(n) \left( \int_M |H|^2 v^2 \, d\mu \right)^{\frac{n-2}{2(n-1)}}\left( \int_M v^{\frac{2n}{n-2}} \, d\mu \right)^{\frac{n-2}{2(n-1)}}\\
\leq & \,C(n)\, \left( \int_M |\nabla v|^2\, d\mu + \int_M H^2 v^2\, d\mu \right) + \frac{1}{2}\left( \int_M v^{\frac{2n}{n-2}} \, d\mu\right)^{\frac{n-2}{n}},
\end{align*}
where in the last inequality we have used the Young's inequality
$$ ab \leq \delta a^p + \delta^{-q/p} b^q$$
for any $a,b,\delta>0$ and $p,q>1$ with $\frac{1}{p} + \frac{1}{q} = 1$. This yields \eqref{eq-MS-1}.

\item[(ii)]
First apply H\"older inequality,  then Michael-Simon's inequality \eqref{eq-MS}  to $w = v^2$ to get
\begin{eqnarray*}
\int_M v^2\, d\mu &\le& C\,\left(\int_M v^4\, d\mu\right)^{\frac 12}
\le C\left( \int_M|\nabla v| |v|\, d\mu + \int_M H v^2\, d\mu\right) \\
&\le& C\left(\, \left(\int_M |\nabla v|^2\, d\mu\right)^{\frac 12} + \left(\int_M H^2 v^2\, d\mu\right)^{\frac 12}\right)\, \left(\int_M v^2\, d\mu\right)^{\frac 12}.
\end{eqnarray*}
This finally implies
$$\int_M v^2 \, d\mu \le C\, \left(\int_M |\nabla v|^2\, d\mu + \int_M v^2H^2\, d\mu\right).$$
\end{enumerate}
\ep

\begin{prop} \label{key1}
Let $M_t^n \subset \mathbb{R}^{n+1}, n\geq 2$ be a smooth compact solution to the {\mcf} \eqref{eq-mcf} for $t\in [0,T)$
with $T < \infty$. Assume that
\be\label{Condition-1}
\max \left\{ \max_{M_0}|A|\,, \,\int_{M_0}|\nabla^m A|^2\,d\mu \right\}\leq  \Lambda_0
\ene
for some $\Lambda_0 \gg 1$ and all $m \in [1, \hat{m}]$ for some fixed $\hat{m} \gg 1$.  Then there exist an
$\epsilon = \epsilon(n, |M_0|, \Lambda_0)>0$, $T_1 = T_1(\Lambda_0) \in (0,1)$, $C_1 = C_1 (n, |M_0|, \Lambda_0)$ and some universal constant $\alpha \in (0,1)$ such that if
\be\label{Condition-2}
\int_{M_0}|\text{\AA}|^2 \,d\mu < \epsilon\,,
\ene
then for all $t\in [0,T_1]$ we have
\be\label{estimate-1}
\max_{M_t}|A| \leq  2\Lambda_0
\ene
and
\be\label{estimate-2}
\max_{M_t}|\text{\AA}| \leq C_1 \epsilon^{\alpha}\,.
\ene
\end{prop}

\bp
By (ii) of Corollary \ref{Evol-A-H} we have
\beq
\frac{\partial}{\partial t} |A|^2 \leq  \Delta |A|^2 + 2|A|^4 \quad \text{on }\, M_t \text{ for all } t\in [0,T)\,.
\eeq
Then by the {\maxp} (Theorem \ref{maxPrin}), we have:
\beq
\max_{M_t}|A| \leq \frac{1}{\sqrt{-2t  + \Lambda_0^{-2}}} \text{ for all } t\in [0,T)\,.
\eeq
Choose $T_1 \leq \frac{3}{8\Lambda_0^2}  \ll 1 $ so that $\max_{M_t}|A| \leq 2 \Lambda_0$ for all $t\in [0,T_1]$. Integrating  equation \eqref{hdA-2} over $M_t$, and using Hamilton's interpolation inequality for tensors (Theorem \ref{Hamiltonlemma2}), we have the
{\ub} on all  higher order derivatives of $A$, which only depends on $n$ and $\Lambda_0$ (more precisely on
$\max_{M_t}|A|$ and the initial bound on the $L^2$-norms of all the derivatives of $A$ in \eqref{Condition-1}).
In particular, for all $m \in [1,\hat{m}]$, we have:
\be\label{gradientbound1}
\max_{M_t}|\nabla^m A|  \leq C(n, \Lambda_0) \quad \text{for }\, t\in [0,T_1]\,,
\ene
c.f. \cite[Lemma 8.3]{Hui84}.

Now we integrate the {\ee} for $|\text{\AA}|^2$, namely, the equation (iii) of  Corollary \ref{Evol-A-H} over $M_t$ for $t\in [0,T_1]$, to get
$$
\frac{\partial}{\partial t} \int_{M_t}|\text{\AA}|^2 \, d\mu  + \int_{M_t}|\text{\AA}|^2 H^2 \, d\mu = -  2 \int_{M_t} |\nabla \text{\AA}|^2 \, d\mu + 2\int_{M_t} |A|^2 |\text{\AA}|^2 \, d\mu\,,
$$
and therefore
\be\label{good-1}
\frac{\partial}{\partial t} \int_{M_t}|\text{\AA}|^2 \, d\mu \leq 8 \Lambda_0^2 \int_{M_t}|\text{\AA}|^2 \, d\mu \quad \text{for all } t\in [0,T_1]\,.
\ene
Using \eqref{good-1} and the assumption that $\int_{M_0}|\text{\AA}|^2 \,d\mu < \epsilon$ yield to
\be\label{good2}
 \int_{M_t}|\text{\AA}|^2 \, d\mu \leq \epsilon e^{8 \Lambda_0^2t} \leq \sqrt{e}\epsilon \leq 2\epsilon \quad \text{for all } t\in [0,T_1]\,.
\ene
By Hamilton's interpolation inequality (Theorem ~\ref{Hamiltonlemma2} for $r=1,p=q=2$) we have
\be\label{good-eq-6}
\int_{M_t} |\nabla \text{\AA}|^2\,d\mu \leq 2 \left(\int_{M_t} |\text{\AA}|^2\,d\mu\right)^{\frac{1}{2}}
\left(\int_{M_t} |\nabla^2 \text{\AA}|^2\,d\mu\right)^{\frac{1}{2}} \leq C(n, \Lambda_0)\epsilon^{\frac{1}{2}}\,,
\ene
where we used $|\nabla^2 \text{\AA}|\leq C(n)|\nabla^2 A|$ and \eqref{gradientbound1}. In fact, using \eqref{gradientbound1} and
applying Theorem \ref{Hamiltonlemma2} part (i) inductively for $r = 1$, $p = q = 2$ and $\Omega := \nabla^{m-1} \text{\AA}$ , we have for all $m \in [1,\hat{m}]$,
\be
\int_{M_t} |\nabla^m \text{\AA}|^2\,d\mu \leq C(n, m, \Lambda_0)\epsilon^{\frac{1}{2^{m-1}}} \quad \text{for all } t\in [0,T_1]\,.
\ene
This together with Theorem \ref{Hamiltonlemma2} part (ii) (for $\Omega := \nabla^{m-1}\text{\AA}$, $i = 1$ and $j = \frac p2$) imply that, for all $t\in [0,T_1]$,
\beq
\int_{M_t} |\nabla^m \text{\AA}|^p\,d\mu \leq C(n, m, p, \Lambda_0)\epsilon^{\frac{1}{2^{m+p/2-2}}} \quad
\text{for all }m \in [1,\hat{m}] \text{ and } p<\infty.
\eeq
By the standard Sobolev embedding theorems (see e.g. \cite[\S 2]{Aub98}) we have that for some universal constant
$\alpha \in (0,1)$, all $m \in [1,\hat{m}-1]$ and $t\in [0,T_1]$,
\be\label{good3}
\max_{M_t}|\nabla^m \text{\AA}| \leq C(n, m,\Lambda_0)\epsilon^{\alpha}.
\ene
Now observe that  $\frac{d}{dt}d\mu = -H^2\, d\mu$, implying that the $|M_t| \leq |M_0|$, where $|M_t|$ is the area of $M_t$. In particular, \eqref{good2}, \eqref{good3} and Lemma \ref{top} (note that $|H|^2\leq n|A|^2 \leq 4n\Lambda_0^2$ in $[0,T_1]$) together yield a bound on $\max_{M_t}|\text{\AA}|$ in terms of $n, |M_0|, \Lambda_0$ and $\epsilon$ for all $t\in [0,T_1]$, that is, \eqref{estimate-2}.
\ep

From now on we will use the same symbol $\epsilon$ for a small constant and symbol $C$ for a uniform constant depending only on $n$, $|M_0|$ and $\Lambda_0$ (possibly $c_0$). We have the immediate consequence of Proposition \ref{key1}, the uniform bound on  mean curvature along the flow and the fact that $|A|^2 = |\text{\AA}|^2 + \frac{1}{n} H^2$.

\begin{cor}
\label{cor-bound-A}
Let $M_t^n  \subset \mathbb{R}^{n+1}, n\geq 2$ be a smooth compact solution to the {\mcf} \eqref{eq-mcf} for $t\in [0,T)$
with $T < \infty$, and $\sup_{M_t} |H| \le c_0$ for all $t\in [0,T)$. Assume that
\be
\max \left\{ \max_{M_0}|A|\,, \,\int_{M_0}|\nabla^m A|^2\,d\mu \right\}\leq  \Lambda_0
\ene
for some $\Lambda_0 \gg c_0$ and all $m \in [1, \hat{m}]$ for some fixed $\hat{m} \gg 1$.  Then there exist an
$\epsilon > 0$ and  $T_1 = T_1(\Lambda_0) \in (0,1)$, such that
$$\max_{M_t} |A| \leq \frac{\Lambda_0}{2} \leq \Lambda_0 \quad \mbox{for all} \,\,\, t\in (0,T_1].$$
\end{cor}

In what follows we want to show that the $\int_{M_t} |\text{\AA}|^2\, d\mu$ stays small along the flow so that we can use iterative type of arguments.

\bl \label{key3}
Let $M_t^n \subset \mathbb{R}^{n+1}, n\geq 2$ be a smooth compact solution to the {\mcf} \eqref{eq-mcf} with $\sup_{M_t} |H| \le c_0$ for all $t\in [0,T)$, where  $T$ is the singular time and  $T_1$ is as in Proposition \ref{key1}. Then there exists an
$\epsilon$ such that if
 \be
\int_{M_0}|\text{\AA}|^2 \,d\mu < \epsilon,
\ene
then
\be \label{MS-est3-1}
\frac{d}{dt} \int_{M_t} |\text{\AA}|^2\, d\mu \leq 0\, \qquad \mbox{for all} \,\,\, t\in [0,T_1]\,.
\ene
\el

\bp
Using \eqref{good3} for $m=1$ and the inequality (see \cite[Lemma 2.2]{Hui84})
$$
|\nabla H|^2 \leq \frac{n(n+2)}{2(n-1)} |\nabla \text{\AA}|^2\,,
$$
we obtain for any $t\in [0,T_1]$
$$
\max_{M_t}|\nabla H| \leq C \epsilon^{\alpha} \,.
$$
Therefore there exists $\eta = C \epsilon^{\frac{\alpha}{2}}>0$ such that if the $\max_{M_{t_0}} |H| \geq \eta$ for some $t_0\in [0,T_1]$ then
$$
\min_{M_{t_0}} |H| \geq \frac{\eta}{2}>0\,,
$$
which follows immediately as an application of Lemma \ref{top} (and we use again that $|H|^2\leq  n|A|^2 \leq 4n\Lambda_0^2$ in $[0,T_1]$ by Proposition \ref{key1}). The smoothness of $H$ implies that $H$ does not change sign on $M_{t_0}$. Since there is no closed hypersurface with strictly negative mean curvature, it follows that $H \geq \frac{\eta}{2}> 0$ on $M_{t_0}$. Now using \eqref{estimate-2} and possibly choosing $\epsilon$ even smaller one sees that the principle curvatures $\kappa_i >0$ and the flow stays strictly convex for all $t \ge t_0$. In this case, using the results from \cite{Hui84} we know the flow contracts into a round point so that $H$ must blow up uniformly everywhere. This contradicts with the assumption $\sup_{M_t} |H| \leq c_0$ for all $t\in [0,T)$. Therefore we have $\max_{M_t} |H| \leq \eta = C \epsilon^{\frac{\alpha}{2}}$ for all $t\in [0,T_1]$. Then using \eqref{estimate-2} we have for all $t\in [0,T_1]$
$$
\max_{M_t}|A| \leq C \epsilon^{\frac{\alpha}{2}}\,.
$$
Combining this with  (iii) of Corollary \ref{Evol-A-H} and integrating over $M_t$ yield to
\begin{eqnarray*}
\frac{d}{dt}\int_{M_t} |\text{\AA}|^2\, d\mu + \int_{M_t} |\text{\AA}|^2 H^2\, d\mu &=& -2\int_{M_t} |\nabla\text{\AA}|^2\, d\mu + 2\int_{M_t} |A|^2 |\text{\AA}|^2\, d\mu \\
&\le& -2\int_{M_t} |\nabla\text{\AA}|^2\, d\mu + 2C \epsilon^{\alpha} \int_{M_t} |\text{\AA}|^2\, d\mu\,.
\end{eqnarray*}

In the case $n = 2$ the proof of \eqref{MS-est3-1} now follows by Lemma \ref{MSineq2} (ii) applied to $v=|\text{\AA}|$ (with Kato's inequality $|\nabla |\text{\AA}||\leq |\nabla \text{\AA}|$) and by choosing  $\epsilon$ small.

If $n > 2$, using H\"older's inequality we obtain for $t\in [0,T_1]$
$$
\frac{d}{dt}\int_{M_t} |\text{\AA}|^2\, d\mu + \int_{M_t} |\text{\AA}|^2 H^2\, d\mu \le -2\int_{M_t} |\nabla\text{\AA}|^2\, d\mu + C\epsilon^{\alpha} \left(\int_{M_t} |\text{\AA}|^{\frac{2n}{n-2}}\, d\mu\right)^{\frac{n-2}{n}}.
$$
This together with Lemma \ref{MSineq2} (i) yields to \eqref{MS-est3-1} by choosing $\epsilon$ sufficiently small.

\ep

Next we prove the $L^2$-norms of higher order derivatives of the second fundamental form $A$ are monotonically non-increasing.

\bl
\label{lem-der}
 Let $M_t^n \subset \mathbb{R}^{n+1}, n\geq 2$ be a smooth compact solution to the {\mcf} \eqref{eq-mcf}  and let $T_1$ and $\hat{m}$ be as in Proposition \ref{key1}. Then there exists an
$\epsilon$ such that if $\int_{M_0}|\text{\AA}|^2 \,d\mu < \epsilon$, then for all $m\in [1,\hat{m}-1]$ we have
\be
\label{MS-est4}
\frac{d}{dt} \int_{M_t} |\nabla^m A|^2\, d\mu \leq 0\,.
\ene
\el
\bp
Using \eqref{hdA-3}, we can finish the proof of Lemma \ref{lem-der} following the proof of Lemma \ref{key3} (recalling Kato's inequality $|\nabla |\nabla^m A|| \leq |\nabla^{m+1} A|$ and replacing $|\text{\AA}|$ by $|\nabla^m A|$, since the evolution equations of the two quantities are almost identical).
\ep

We will end this section by proving Theorem \ref{MainThm}.

\bp[Proof of Theorem \ref{MainThm}] Assume that there exists $\Lambda_0 \gg c_0 >0$ sufficiently large so that condition \eqref{Condition-1} is satisfied. Then by Proposition \ref{key1}, Lemma \ref{key3}, Lemma \ref{lem-der} and Corollary \ref{cor-bound-A} there exist $T_1 =  T_1(\Lambda_0) \in (0,1)$ and  $\epsilon >0$ sufficiently small such that if
$$
\int_{M_0}|\text{\AA}|^2 \,d\mu  <  \epsilon\,,
$$
then only two things can happen:
\begin{enumerate}
\item
the flow becomes strictly mean convex at some time $t_0\in [0,T_1]$. In this case the flow will stay strictly mean convex for all $t\in [t_0, T)$.  By the pinching estimate (see e.g. \cite{HS99}, \cite{Smo98}) we have uniform constant $C > 0$ so that for all $t\in [t_0, T)$
$$|A|^2 \le C( H^2 + 1) \le C( c_0^2 + 1),$$
implying that $T$ is not the singular time and the flow can be smoothly extended past time $T$.
\item
$
\max_{M_t}|A| \leq C \epsilon^{\frac{\alpha}{2}} \le \frac{\Lambda_0}{2} \leq \Lambda_0$,
$$
\frac{d}{dt} \int_{M_t} |\text{\AA}|^2\, d\mu \leq 0 \quad \text{and}\quad \frac{d}{dt} \int_{M_t} |\nabla^m A|^2\, d\mu \leq 0
$$
for all $t\in [0,T_1]$ and $m \in [1,\hat{m}-1]$. These  mean that at time $t= T_1$ we still have  conditions \eqref{Condition-1} and \eqref{Condition-2} holding, i.e.,
$$
\max \left\{ \max_{M_{T_1}}|A|\,, \,\int_{M_{T_1}}|\nabla^m A|^2\,d\mu \right\}\leq  \Lambda_0 \,\,\, \text{and}\,\,\, \int_{M_{T_1}}|\text{\AA}|^2 \,d\mu < \epsilon\,.
$$
Now we can iterate the arguments using Proposition \ref{key1}, Lemma \ref{key3} and Lemma \ref{lem-der} for time intervals of size $T_1 = T_1(\Lambda_0) >0$. Either at some time $t_0$ the flow becomes strictly mean convex in which case we apply the reasoning in (1) and we are done or, since $T_1$ is of a uniform size, after finitely many iterations we reach time $T$ and obtain
$$
\max_{M_t}|A| \leq \Lambda_0 \quad \text{at time } t=T\,.
$$
This means the flow can be smoothly extended past $T$. \qedhere
\end{enumerate}
\ep

\section{The normalized mean curvature flow}
\label{sec-conv}

In the previous section we show that if  the mean curvature flow starts from a hypersurface with small traceless second fundamental form, the mean curvature has to blow up at a singular time. An interesting question to ask is  whether this flow becomes strictly convex at some time before the singularity occurs and therefore becomes extinct in an asymptotically spherical manner as it was shown in \cite{Hui84} (meaning that the normalized mean curvature flow converges as time approaches infinity to a sphere). We  will answer this question by looking at its normalized flow. In fact, in this section we show that such a flow indeed  becomes strictly convex at some time. We shall remark that for the one dimensional curve shortening flow, Grayson proved in \cite{G87} that the flow starting from arbitrary embedded closed curve in $\mathbb{R}^2$ becomes strictly convex at some time and therefore by Gage-Hamilton's results in \cite{GH86} the flow shrinks to a round point. Without loss of generality we assume that the origin is always in the region enclosed by the evolving hypersurfaces for all times $0\leq t <T$. We are going to normalize the mean curvature flow by keeping the total area of the hypersurface $M_t$ constant, as in \cite{Hui84}. Namely, we multiply the solution $F$ of \eqref{eq-mcf} at each time $0\leq t < T$ with a positive constant $\psi(t)$ such that the total area of the hypersurface $\tilde{M}_t$ given by
$$
\tilde{F}(\cdot, t) = \psi(t)\cdot F(\cdot, t)
$$
is equal to the total area of $M_0$. Define $\tilde{t}(t) = \int_0^t \psi(\tau)d\tau$. One can derive the normalized evolution equation for $\tilde{F}$ on a different maximal time interval $0\leq \tilde{t} < \tilde{T}$:
\be \label{5}
\frac{\partial \tilde{F}}{\partial \tilde{t}} = - \tilde{H} \tilde{\nu} + \frac{1}{n} \tilde{h}\tilde{F}\,,
\ene
where $\tilde{h} = \int_{\tilde{M}_{\tilde{t}}} \tilde{H}^2 d\tilde{\mu}\,/\, \int_{\tilde{M}_{\tilde{t}}}  d\tilde{\mu}$. Since $|\tilde{\text{\AA}}|^2 = \psi^{-2}|\text{\AA}|^2$, by Lemma 9.1 of \cite{Hui84}, we have
\be \label{6}
\frac{\partial}{\partial \tilde{t}} |\tilde{\text{\AA}}|^2 = \Delta |\tilde{\text{\AA}}|^2 - 2|\nabla \tilde{\text{\AA}}|^2 + 2|\tilde{A}|^2|\tilde{\text{\AA}}|^2 -\frac{2}{n} \tilde{h} |\tilde{\text{\AA}}|^2\,.
\ene
Now since
$$
\frac{\partial}{\partial \tilde{t}} d \tilde{\mu} = (\tilde{h} - \tilde{H}^2) d\tilde{\mu}\,,
$$
integrating \eqref{6} over $\tilde{M}_{\tilde{t}}$ we have
\begin{align}\label{8}
&\frac{\partial}{\partial  \tilde{t}} \int_{\tilde{M}_{\tilde{t}}} |\tilde{\text{\AA}}|^2 d\tilde{\mu} + \int_{\tilde{M}_{\tilde{t}}}  \tilde{H}^2 |\tilde{\text{\AA}}|^2 d\tilde{\mu} \notag\\
= & -2 \int_{\tilde{M}_{\tilde{t}}} |\nabla \tilde{\text{\AA}}|^2 d\tilde{\mu} + 2 \int_{\tilde{M}_{\tilde{t}}} |\tilde{A}|^2|\tilde{\text{\AA}}|^2 d\tilde{\mu} + \frac{n-2}{n} \tilde{h} \int_{\tilde{M}_{\tilde{t}}} |\tilde{\text{\AA}}|^2 d\tilde{\mu} \,.
\end{align}

Next we prove  Theorem \ref{thm-conv} which claims the convergence of the normalized mean curvature flow to the sphere if the initial hypersurface is in $L^2$ sense close to the sphere.

\bp[Proof of Theorem \ref{thm-conv}]
We argue similarly  as in the proof of Theorem \ref{MainThm}. First note that the area $|M_t|$ is non-increasing along the mean curvature flow \eqref{eq-mcf}, and thus the normalizing factor $\psi(t)\geq 1$ for all $t\in[0,T)$. Assume $\Lambda_0$ is such that
$$\max\left\{\max_{\tilde{M}_0} |\tilde{A}|, \,\int_{\tilde{M}_0} |\nabla^m \tilde{A}|^2\, d\tilde{\mu} \right\} \le \Lambda_0,$$
for some $\Lambda_0 \gg 1$ and all $m \in [1,\hat{m}]$ for some fixed $\hat{m} \gg 1$.
As in the proof  of Proposition \ref{key1}, using the evolution for $|\tilde{A}|^2$ (which is the same as \eqref{6} after replacing $|\tilde{\text{\AA}}|^2$ by $|\tilde{A}|^2$) we have that $\max_{\tilde{M}_{\tilde{t}}} |\tilde{A}| \le 2\Lambda_0$ for all $\tilde{t}\in [0,\tilde{T_1}]$ where $\tilde{T_1} = \tilde{T_1}(\Lambda_0) \in (0,1)$. We have two possible scenarios.
\begin{enumerate}
\item
There exists an $\eta > 0$ to be determined later so that the $\max_{\tilde{M}_{\tilde{t}_0}} |\tilde{H}| \ge \eta$ for some $\tilde{t}_0 \in [0,\tilde{T_1}]$.
\item
For all $\tilde{t}\in [0,\tilde{T_1}]$, we have $\max_{\tilde{M}_{\tilde{t}}} |\tilde{H}| \le \eta$.
\end{enumerate}
In the first case (1), since $\max_{\tilde{M}_{\tilde{t}}} |\tilde{H}| \le 2\sqrt{n} \Lambda_0$ for all $\tilde{t}_0\in [0,\tilde{T_1}]$ (to apply Lemma \ref{top}) the same proof as in Proposition \ref{key1}  yields to
\begin{equation}
\label{eq-small-tilde}
\max_{\tilde{M}_{\tilde{t}}} |\nabla^m \tilde{\text{\AA}}| \le C(n,m,\Lambda_0) \epsilon^{\alpha},
\end{equation}
for all $\tilde{t} \in [0,\tilde{T_1}]$, all $m \in [0,\hat{m}]$ where $\alpha \in (0,1)$ is some universal constant and $\epsilon$ is sufficiently small. Since $|\nabla\tilde{H}| \le \left(\frac{n(n+2)}{2(n-1)}\right)^{1/2} |\nabla\tilde{\text{\AA}}| \le C(n, \Lambda_0) \epsilon^{\alpha}$ and $\max_{\tilde{M}_{\tilde{t}_0}} |\tilde{H}| \ge \eta$, if we choose $\eta = c_2(n, |M_0|, \Lambda_0) \epsilon^{\frac{\alpha}{2}}$, we get
$$\min_{\tilde{M}_{\tilde{t}_0}} \tilde{H} \ge \frac{\eta}{2},$$
by the same arguments as in Proposition \ref{key1}. Combining this with \eqref{eq-small-tilde} for $m = 0$, if $\epsilon > 0$ is chosen sufficiently small, yields to strict convexity of $\tilde{M}_{\tilde{t}_0}$. One can now apply the results in \cite{Hui84} to get the long time existence of the normalized mean curvature flow and its exponential convergence to a round sphere.

\noindent In the second case (2), the same arguments as in the proof of Proposition \ref{key1}  yield to
\be \label{9}
\max_{\tilde{M}_{\tilde{t}}}|\tilde{A}| = \psi^{-1}(\tilde{t})\max_{M_{\tilde{t}}}|A| \leq \max_{M_{\tilde{t}}}|A| \leq C \epsilon^{\frac{\alpha}{2}} \le \frac{\Lambda_0}{2} \leq  \Lambda_0
\ene
for all $\tilde{t}\in [0,\tilde{T_1}]$. This also implies
\be \label{11}
\max_{\tilde{t}\in [0,\tilde{T_1}]}|\tilde{h}(\tilde{t})| \leq C \epsilon^{\alpha}\,.
\ene
 Moreover, combining \eqref{8}, \eqref{9} and \eqref{11} and following the proof of Lemma \ref{key3} we get
\be\label{10}
\frac{\partial}{\partial  \tilde{t}} \int_{\tilde{M}_{\tilde{t}}} |\tilde{\text{\AA}}|^2 d\tilde{\mu} \leq - \int_{\tilde{M}_{\tilde{t}}} |\nabla \tilde{\text{\AA}}|^2 d\tilde{\mu} - \frac{1}{2} \int_{\tilde{M}_{\tilde{t}}} \tilde{H}^2 |\tilde{\text{\AA}}|^2 d\tilde{\mu} \leq 0
\ene
for all $\tilde{t}\in [0,\tilde{T_1}]$ if we choose $\epsilon$ sufficiently small. Similarly it can be shown that in this case we also have
$$
\frac{\partial}{\partial  \tilde{t}} \int_{\tilde{M}_{\tilde{t}}} |\tilde{\nabla}^m \tilde{\text{\AA}}|^2 d\tilde{\mu} \leq 0
$$
for all $\tilde{t}\in [0,\tilde{T_1}]$ and $m \in [1,\hat{m}-1]$. In particular the analysis in case (2) implies that
$$\max\left\{ \max_{\tilde{M}_{\tilde{T_1}}} |\tilde{A}|, \, \int_{\tilde{M}_{\tilde{T_1}}} |\nabla^m \tilde{A}|^2\, d\mu\right\} \leq \Lambda_0 \quad \mbox{and} \quad \int_{\tilde{M}_{\tilde{T_1}}}|\tilde{\text{\AA}}|^2\, d\mu < \epsilon\,.$$

We can then iterate the arguments for another uniform size $\tilde{T_1}>0$ for infinitely many times and one sees that the flow exists for all time $\tilde{t}>0$. Moreover, $\tilde{A}$ and all its derivatives are uniformly bounded along the normalized flow. If for some $\tilde{t}_0$ we happen to be in case (1) we are done. Otherwise assume case (2) happens for all $\tilde{t} \in [0,\infty)$. In this case there is a uniform constant $C$ depending only on $M_0$ and $n$ so that $|\nabla^m \tilde{A}|(\cdot,\tilde{t}) \le C$ for all $m \in [0,\hat{m}]$ and all $\tilde{t} \in [0,\infty)$. This implies for any sequence $\tilde{t}_i\to\infty$ there exists a subsequence so that the normalized mean curvature flow converges in the $C^k$-topology to a smooth closed hypersurface $\tilde{M}_{\infty}$. The monotonicity formula \eqref{10} forces
$$|\nabla\tilde{\text{\AA}}_{\infty}| \equiv 0 \qquad \mbox{and}\qquad  |\tilde{\text{\AA}}_{\infty}| |\tilde{H}_{\infty}| \equiv 0,$$
on $\tilde{M}_{\infty}$.  The first equation implies $\tilde{\text{\AA}}_{\infty} \equiv const$. If  $const = 0$, i.e., $\tilde{M}_{\infty}$ is a closed umbilic hypersurface in $\mathbb{R}^{n+1}$, then it is well-known that $\tilde{M}_{\infty}$ must be an embedded round sphere, see e.g. \cite{Ger}. If $const \neq 0$, the second equation forces $\tilde{H}_{\infty} \equiv 0$, which is impossible since there are no closed minimal hypersurfaces in $\mathbb{R}^{n+1}$. Therefore there exists a sufficiently large  time at which the flow becomes strictly convex and then converges exponentially to a round sphere with area $|M_0|$ by the results in \cite{Hui84}.
\ep

\bibliographystyle{amsalpha}
\bibliography{ref-2dmcf}
 \end{document}